# Some recent advances on the RBF

W. Chen

*Scientific Computing Department, Simula Research Lab., P. O. Box. 134, 1325 Lysaker, Norway (wenc@simula.no)*

## Abstract

This paper aims to survey our recent work relating to the radial basis function (RBF) and its applications to numerical PDEs. We introduced the kernel RBF involving general pre-wavelets and scale-orthogonal wavelets RBF. A dimension-independent RBF error bound was also conjectured. The centrosymmetric structure of RBF interpolation matrix under symmetric sample knots was pointed out. On the other hand, we introduced the boundary knot method via nonsingular general solution and dual reciprocity principle and the boundary particle method via multiple reciprocity principle. By using the Green integral we developed a domain-type Hermite RBF scheme called the modified Kansa method, which significantly reduces calculation errors around boundary. To circumvent the Gibbs phenomenon, the least square RBF collocation scheme was presented. All above discretization schemes are meshfree, symmetric, spectral convergent, integration-free and mathematically very simple. The numerical validations are also briefly presented.

## 1. Introduction

Since Kansa's pioneer work [1] in 1990, the research on the RBF for numerical PDEs has become very active. Nevertheless, Nardini and Brebbia [2] have actually much earlier applied the RBF concept to develop the dual reciprocity BEM (DR-BEM) without using the term "RBF". Unlike most of other meshfree schemes, the RBF approach doe not require using the moving least square (MLS). The RBF invariably involves only one-dimensional distance and is therefore inherently meshfree and independent of geometric complexity.

   This paper aims to survey our latest progresses on the RBF and related applications to numerical PDEs. For complete references see the respective author's paper. For other new advances, especially for so-called fast RBF, see [3]. In section 2, the kernel RBF [4-9] is discussed, which creates operator-dependent RBFs via the fundamental and general solutions. The popular thin plate spline (TPS) and multiquadratics (MQ) are special cases of the kernel RBF. It is noted that the distribution theory and integral equation are important to research the RBF. In section 3, the newly-discovered high-order fundamental and general solutions of convection-diffusion, Winker plate and Burger plate are also given [10-12]. Section 4 is concerned with the PDE-based RBF scale-orthogonal wavelets and general pre-wavelets [5,6,13,14] using the fundamental and general solutions of differential operators, which combine the strengths of both the RBF and wavelets. In section 5 we discuss the conjecture error estimate [15] and centrosymmetric matrix structure under symmetric sample knots of the RBF interpolation [6]. Finally, section 6 introduces a few new RBF discretization schemes of boundary and domain types [9,20], which include the boundary knot method (BKM) [4-6,9-11,16-20], boundary particle method (BPM) [9,11,12,20,21], modified Kansa method (MKM) [9,20], and least square RBF collocation method (LSRCM) [9]. Among them, the BKM applies the nonsingular general solution and dual reciprocity principle, while the BPM uses the multiple reciprocity principle. The MKM significantly reduces the errors of the Kansa's method at boundary-adjacent region. The LSRCM circumvents the Gibbs phenomenon via the Least square approximation instead of the interpolation.

## 2. Kernel RBF

The origins of the traditional RBFs (except the TPS for 2D biharmonic equation) have little to do with certain PDE problem. Based on the second Green identity, Chen [4-9] presented the kernel RBF-creating strategy. Let us consider the following example without loss of generality:

$$\Re\{u\} = f(x), \qquad x \in \Omega, \qquad (1)$$

$$u(x) = R(x), x \subset S_u \text{ and } \frac{\partial u(x)}{\partial n} = N(x), \quad x \subset S_T, \quad (2a,b)$$

where $x$ means multi-dimensional independent variable, and $n$ is the unit outward normal. The Green integral solution of Eqs. (1,2a,b) is given by

$$u(x) = \int_\Omega f(z) u^*(x,z) d\Omega + \int_\Gamma \left\{ u \frac{\partial u^*(x,z)}{\partial n} - \frac{\partial u}{\partial n} u^*(x,z) \right\} d\Gamma, \qquad (3)$$

where $u^*$ is the fundamental solution of differential operator $\Re\{\}$. $z$ denotes source point. The numerical analog of Eq. (3) is stated as

$$u(x) \cong \sum_{k=1}^{N+L} \omega(x,x_k) f(x_k) u^*(x,x_k) + \sum_{k=N+1}^{N+L} \varrho(x,x_k) \left[ u \frac{\partial u^*(x,x_k)}{\partial n} - \frac{\partial u}{\partial n} u^*(x,x_k) \right], \qquad (4)$$

where $\omega$ and $Q$ are the weight functions dependent on the integral schemes, $N$ and $L$ respectively denote the numbers of domain and boundary knots. It is observed that (4) is a RBF solution if $\omega$ and $Q$ are only dependent on the distance variable. Five types of the kernel RBF were proposed. The first is to apply $r^{2m}$ augmented term to enhance the smoothness and ensures sufficient degree of differential continuity [4-9]. The TPS is a notable example in this regard. The second strategy is simply the higher-order fundamental and general solutions [9,20]. The third approach is to replace distance variable $r$ in fundamental or general solutions by $\sqrt{r^2+c^2}$ [4-6,9], where $c$ is the shape parameter. Following the basic idea of the corrected reproducing kernel approximation [22], the forth is to construct the compactly-supported kernel RBF [9]. For instance, the MQ can be used as a correction function to decide local optimal shape parameter by establishing the reproducing conditions. The fifth kernel RBF is the time-space RBF [4-7,9] which employs transient fundamental solution and general solution, e.g. the characteristic RBF $\varphi(r-ct)$ for hyperbolic wave problems.

The operator-dependent kernel RBF is strongly recommended in numerical PDE and data processing. For example, $e^{\pm\alpha r}$ [9] is much better than Guassian, TPS or MQ for diffusion and convection-diffusion problems due to its underlying approximation to their fundamental and general solutions.

## 3. New high-order fundamental and general solutions

Table 1 displays the newly-found high-order fundamental and general solutions of convection-diffusion [11,12], vibration, Winkler and Burger plates [10]. These solutions are verified via the computer software "Maple".

Table 1. $M$-order general solutions, where $n=2,3$ denotes dimensionality.

|  | Operators | General solutions ($m=0,1,2…$) |
|---|---|---|
| Convection-diffusion | $D\nabla^2 u - v\bullet\nabla u - \kappa u$ | $u_m^\#(r) = Q_m(\mu r)^{-n/2+1+m} e^{\frac{v\cdot r}{2D}} I_{n/2-1+m}(\mu r)$ |
| Vibration plate | $\nabla^4 u - \lambda^2 u$ | $u_m^\#(r) = (r\sqrt{\lambda})^{-n/2+1+m}\left(A_m J_{n/2-1+m}(\sqrt{\lambda}r)\right.$ $\left. + B_m I_{n/2-1+m}(\sqrt{\lambda}r)\right)$ |
| Winkler plate | $\nabla^4 u + \kappa^2 u$ | $u_m^\#(r) = (r\sqrt{\kappa})^{-n/2+1+m}(C_m ber_{n/2}(r\sqrt{\kappa})$ $+ D_m bei_{n/2}(r\sqrt{\kappa}))$, $m$ odd <br> $u_m^\#(r) = (r\sqrt{\kappa})^{-n/2+1+m}(C_m ber_{n/2-1}(r\sqrt{\kappa})$ $+ D_m ber_{n/2-1}(r\sqrt{\kappa}))$, $m$ even |
| Burger plate | $\nabla^4 u - \mu^2\nabla^2 u$ | $u_m^\#(r) = E_m r^{2m-2} + F_m(\mu r)^{-n/2+1+m} I_{n/2-1+m}(\mu r)$ |

The general solution satisfies the differential equation, but contrast to the fundamental solution, at origin is equal to a limited value rather than zero and infinity. For the convection-diffusion operator, $I$ represents the modified Bessel function of the second kind and $Q_m=Q_{m-1}/(2*m*\mu^2)$, $Q_0=1$. The high-order fundamental solutions are obtained by replacing $I$ by the first kind modified Bessel function.

$$\mu = \left[\left(|v|/2D\right)^2 + \kappa/D\right]^{\frac{1}{2}}. \tag{5}$$

For vibration, Winkler and Burger plates, ber and bei respectively represent the Kelvin and modified Kelvin functions of the first kind. $A_m$ to $F_m$ are constant coefficients which will be detailed in a subsequent paper. It is worth pointing out that the formulas given in Table 1 for the zero order general solution of Winkler operator is effective for up to five dimensions. The same relations hold with ber, bei replaced by the Kelvin functions of the second kind ker, kei, respectively, for fundamental solutions. The higher order fundamental solution of Burger equation is alternations of the first and second terms of corresponding general solution by the higher order fundamental solutions of Laplace and Helmholtz operators.

## 4. Error estimate and centrosymmetric structures

The existing RBF error bounds do not consider dimension effect. It is very interesting to observe that there exist the same error behaviors between some RBFs and quasi-Monte Carlo (QMC) or Monte Carlo method (MCM). For example, error bounds for the linear RBF and the classical Monte Carlo method are the same $O(M^{-1/2})$ which $M$ is the number of nodes, while error bounds for the TPS and QMC in 2D problems are the same $O(M^{-1}(\log M))$. In fact, the RBF and QMC or MCM have some close relationship on the grounds of numerical integration [15]. By analogy with the error estimate of the QMC and MCM, we intuitively proposed error bound conjecture of the RBF [15]: $O(M^{-\eta}(\log M)^{d-1})$, where $d$ is dimensionality. The notorious dimension curse in other numerical techniques can be characterized by $err = O\left(M^{-\kappa/d}\right)$. The RBF has visible advantages in accuracy for higher dimensional problems.

If the nodes are symmetrically sampled, the RBF matrix is either centrosymmetric for even order derivative or skew-centrosymmetric for odd order derivative, [6,20], which are defined as

$$r_{ij} = r_{N+1-i,N+1-j} \quad \text{and} \quad r_{ij} = -r_{N+1-i,N+1-j}. \tag{6a,b}$$

Centrosymmetric matrices can easily be decomposed into two half-sized matrices. Such factorization leads to a considerable reduction in computing effort for determinant, inversion and eigenvalues (for details see refs in [6]).

## 5. RBF general pre-wavelets and scale-orthogonal wavelets

The RBF is well known for its striking effectiveness in multivariate scattered data approximation. However, in general the RBFs available now lack critical multiscale analysis capability. To handle high-dimensional multiscale analysis, the RBF wavelets are mostly wanted to combine the strengths of both. In last decade much effort has been devoted to non-orthogonal pre-wavelets RBF theory by using some constructive approximation strategies.

### 5.1. RBF general pre-wavelets

Buhmann [23] shows that the MQ is pre-wavelets, where the shape parameter $c$ is seen as the scaling coefficient. In the preceding section 2, the kernel RBF which replaces distance $r$ in fundamental or general solutions by $\sqrt{r^2+c^2}$ should also be understood general pre-wavelets. For instance, numerical experiments with pre-wavelet TPS $r_j^{2m}\ln\sqrt{r_j^2+c_j^2}$ or $\left(r_j^2+c_j^2\right)^m\ln\sqrt{r_j^2+c_j^2}$ manifest spectral convergence as in the MQ [5.6].

It is noted that the fundamental solution used in the BEM involves only the essential part of a complete fundamental solution [24]. The complementary term is often regarded as the nonsingular general solution in terms of the BKM. From this point of view, the shape parameter $c$ can be interpreted as the scaling parameter in the simplified form of the complete fundamental solutions and leads to infinite smoothness at the cost of one complementary term. For instance, the MQ is related with general fundamental solutions of the Laplacian, where the complementary term is a constant. The tricky choose of the shape parameter coincides with the skillful implementation of general fundamental solution. This may also reveal some subtle relationship between the wavelets and physical field theory.

### 5.2. Scale-orthogonal wavelets RBF

By comparing to Fourier series and transform, Chen [13,14] developed the scale-orthogonal RBF wavelet series and transform via the fundamental and general solutions of some typical PDEs. This work, however, is more conjecture and speculation than a complete theory, notably lacking rigorous mathematical analysis. Ref. 14 presented the RBF wavelet analytical solution of transient scalar wave of arbitrary dimensionality and geometry, which, to the best of the author's knowledge, had not been achieved before. It is promising to construct a complete RBF wavelets theory by combing PDE fundamental and general solutions and the central basis function approach.

## 6. RBF numerical discretizations

The Kansa's method [1] is the very first domain-type RBF collocation scheme with easy-to-use merit, but the method lacks symmetric interpolation matrix due to boundary collocation. The Hermite RBF collocation method kills the unsymmetrical drawback. Like the Kansa's method, however, the method suffers relatively lower accuracy in boundary-adjacent region. The method of fundamental solution (MFS) [25], also known as the regular BEM, is a simple and efficient boundary-type RBF scheme, but the controversial artificial boundary outside physical domain hinders its practical applications and causes the unsymmetric interpolation matrix. The global RBF interpolation is also ill-conditioning and susceptible to Gibbs phenomenon amid weak continuity of physical solution. This section introduces a few recent RBF discretization schemes of boundary and domain types [9,20] to overcome the aforementioned shortcomings, all of which are meshfree, symmetric, integration free, and mathematically simple. In addition, both indirect expansion coefficients and direct physical variables can be applied as basic variables within these methods.

### 6.1. Boundary knot method

The BKM is a two-step technique based on the fact that the solution of Eq. (1) can be split into homogeneous solution $u_h$ and particular solution $u_p$. Like the DR-BEM and MFS, the particular solution is evaluated by the RBF and dual reciprocity principle. Unlike them, the homogeneous solution is approximated by nonsingular general solution. For symmetric interpolation,

$$u_h(x) = \sum_{s=1}^{L_D} \lambda_s u^\#(r_s) - \sum_{s=L_D+1}^{L_D+L_N} \lambda_s \frac{\partial u^\#(r_s)}{\partial n}. \tag{7}$$

The BKM numerical formulation is given by

$$\sum_{s=1}^{L_D} \lambda_s u^\#(r_{is}) - \sum_{s=L_D+1}^{L_D+L_N} \lambda_s \frac{\partial u^\#(r_{is})}{\partial n} = R(x_i) - u_p(x_i), \tag{8}$$

$$\sum_{s=1}^{L_D} \lambda_s \frac{\partial u^\#(r_{js})}{\partial n} - \sum_{s=L_D+1}^{L_D+L_N} \lambda_s \frac{\partial^2 u^\#(r_{js})}{\partial n^2} = N(x_j) - \frac{\partial u_p(x_j)}{\partial n}, \tag{9}$$

$$\sum_{s=1}^{L_D} \lambda_s u^\#(r_{ls}) - \sum_{s=L_D+1}^{L_D+L_N} \lambda_s \frac{\partial u^\#(r_{ls})}{\partial n} = u_l - u_p(x_l), \tag{10}$$

where $i$, $j$, and $l$ indicate response knots respectively located on boundary $S_u$, $S_\Gamma$, and domain $\Omega$. Note that inner nodes are usually necessary to ensure accuracy, stability and convergence.

### 6.2. Boundary particle method

The multiple dual reciprocity BEM applies the multiple reciprocity principle to avoid the domain integral without using any inner nodes [26]. The drawback is uneasily used to nonlinear problems and requires higher computing effort. Based on the multiple reciprocity principle and RBF, we developed the meshfree BPM to overcome such difficulties while keeping truly boundary-only merit. The BPM formulation of Eqs. (1,2a,b) is given by

$$\sum_{s=1}^{L_D} \beta_s u_0^\#(r_{is}) - \sum_{s=L_D+1}^{L_D+L_N} \beta_s \frac{\partial u_0^\#(r_{is})}{\partial n} = R(x_i) - u_p^0(x_i)$$
$$\sum_{s=1}^{L_D} \beta_s \frac{\partial u_0^\#(r_{js})}{\partial n} - \sum_{s=L_D+1}^{L_D+L_N} \beta_s \frac{\partial^2 u_0^\#(r_{js})}{\partial n^2} = N(x_j) - \frac{\partial u_p^0(x_j)}{\partial n} \quad (11)$$

$$\sum_{s=1}^{L_D} \beta_s \Re^{n-1}\{u_n^\#(r_{is})\} - \sum_{s=L_D+1}^{L_D+L_N} \beta_s \frac{\partial \Re^{n-1}\{u_n^\#(r_{is})\}}{\partial n} = \Re^{n-2}\{f(x_i)\} - \Re^{n-1}\{u_p^n(x_i)\}$$
$$\sum_{s=1}^{L_D} \beta_s \frac{\partial \Re^{n-1}\{u_n^\#(r_{js})\}}{\partial n} - \sum_{s=L_D+1}^{L_D+L_N} \beta_s \frac{\partial^2 \Re^{n-1}\{u_n^\#(r_{js})\}}{\partial n^2} = \frac{\partial(\Re^{n-2}\{f(x_j)\} - \Re^{n-1}\{u_p^n(x_j)\})}{\partial n}$$
$$n=1,2,\ldots, \quad (12)$$

The practical solution procedure is a reversal recursive process:

$$\beta_k^M \to \beta_k^{M-1} \to \cdots \to \beta_k^0. \quad (13)$$

Since all successive equations have the same interpolation matrices, the LU decomposition algorithm is suitable. The solution at any node is then given by

$$u(x_p) = \sum_{l=0}^{M} \left( \sum_{s=1}^{L_D} \beta_s^l u_m^\#(r_{ps}) - \sum_{s=L_D+1}^{L_D+L_N} \beta_s^l \frac{\partial u_m^\#(r_{ps})}{\partial n} \right). \quad (14)$$

### 6.3. Modified Kansa method

Considering numerical discretization (4) of Green integral solution of Eqs. (1,2a,b), we can construct the following RBF interpolation formula

$$u(x) = \sum_{k=1}^{N+L} \alpha_k L^*\{\varphi(r_k)\} + \sum_{k=1}^{N+L_D} \beta_k \varphi(r_k) + \sum_{k=N+L_D+1}^{N+L} \beta_k \left(-\frac{\partial \varphi(r_k)}{\partial n}\right), \quad (15)$$

where $\Re^*\{\}$ reverses some sign of odd-order derivative in $\Re\{\}$ if the latter is not self-adjoint. Note that the boundary nodes are here interpolated twice. Collocating Eqs. (1,2a,b) via interpolation formula (15) leads to the modified Kansa method which greatly reduces the solution errors around the boundary.

### 6.4. Least square RBF collocation method

For the least square RBF collocation method, the field and source knots are unnecessarily either the same amount or at the same location. The least square approach is applied to the solution of the global and local RBF collocation equations. The merit of the method is to solve the discontinuous problems such as shock. For details see [9].

### 6.5. Numerical experiments to the BKM and BPM

The 2D and 3D irregular geometries tested are illustrated in Figs. 1 and 2 (two-ball cavity: radii=1, center distance=$\sqrt{2}$ ). Except for the specified Neumann boundary conditions shown in Fig. 1 and on the $x$=0 surface in Fig. 2, the otherwise boundary conditions are all of the Dirichlet type. The BKM employs 9 inner nodes for 2D inhomogeneous cases as shown by small crosses in Fig. 1. The $L_2$ norms of relative errors at 364 nodes for 2D and at 500 nodes for 3D are displayed in Tables 2 and 3, respectively.

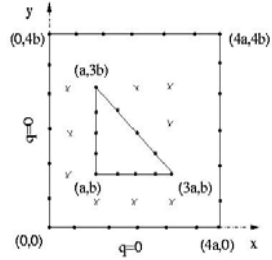 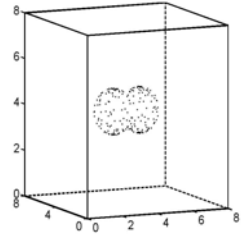

Fig. 1. Configuration of a square with a trigeometric cutout    Fig. 2. A cube with a two-ball cavity

Equations of Helmholtz and convection-diffusion are

$$\nabla^2 u + \gamma^2 u = f(x), \text{ and } D\nabla^2 u(x) - v \bullet \nabla u(x) - \kappa u(x) = g(x), \quad (16a,b)$$

where $v$ denotes velocity vector and $D$ is the diffusivity coefficient, $\kappa$ represents the reaction coefficient. The accurate solutions are

$$u = x^2 \sin x \cos y \quad (17)$$

for 2D inhomogeneous Helmholtz problem and

$$u = \sin x \cos y \cos z \quad \text{and} \quad u = e^{-\sigma x} + e^{-\sigma y} + e^{-\sigma z} \quad (18a,b)$$

for 3D homogeneous Helmholtz and convection-diffusion ($v_x$=$v_y$=$v_z$=-$\sigma$, $\kappa$=0, Pelect number is 24 for $\sigma$=1 and 480 for $\sigma$=20) problems.

Table 2. $L_2$ norm of relative errors for 2D inhomogeneous Helmholtz problems (numbers inside parentheses indicate boundary plus inner nodes).

| BKM (26+9) | BKM (33+9) | BPM (26) | BPM (33) |
|---|---|---|---|
| 1.9e-3 | 9.3e-5 | 2.7e-3 | 6.8e-4 |

Table 3. $L_2$ norm of relative errors for 3D homogeneous Helmholtz and convection-diffusion problems by the BKM (the numbers in parentheses indicate boundary nodes).

| Helmholtz | | Convection-diffusion ($\sigma=1$) | |
|---|---|---|---|
| 4.6e-3 (298) | 1.7e-4 (466) | 9.0e-3 (136) | 2.2e-3 (298) |